\documentclass[12pt]{article}
\usepackage{latexsym}

\usepackage[tbtags]{amsmath}
\usepackage{epsfig,amstext,amssymb,amsthm,latexsym}
\usepackage{amssymb,color}
\usepackage{tikz}
\usetikzlibrary{calc,spy}
\usepackage{pgfplots}
\UseRawInputEncoding

\definecolor{c20}{rgb}{0.,0.7,0.}
\definecolor{c30}{rgb}{0.,0.,1.}
\definecolor{c40}{rgb}{1,0.1,0.7}
\definecolor{c50}{rgb}{1,0,0}

\date{}
\newtheorem{lemma}{Lemma}[section]

\newtheorem{remark}{Remark}[section]

\newtheorem{definition}{Definition}[section]

\makeatletter 
\@addtoreset{equation}{section}
\makeatother 

\begin{document}
\title{A Note on the Borel-Cantelli Lemma}

\author{Narayanaswamy  Balakrishnan \thanks{Department of Mathematics and Statistics, McMaster University,
1280 Main Street West, Hamilton, Ontario, Canada L8S 4K1; email: bala@mcmaster.ca} \and Alexei Stepanov  \thanks{\noindent Higher School of Computer Science and Applied Mathematics, Education and Research Cluster "Institute of High Technology``, Immanuel Kant Baltic Federal University, 
A. Nevskogo  14, Kaliningrad, 236041 Russia; email: alexeistep45@mail.ru}}

\maketitle
\begin{abstract} 
In this short note, we discuss the Barndorff-Nielsen lemma, which is a generalization of well-known Borel-Cantelli lemma. Although the result stated in the Barndorff-Nielsen lemma is correct,  it does not follow from the argument proposed in the corresponding proof. In this note, we show this and  offer  an alternative proof of this lemma. We also propose a new generalization of  Borel-Cantelli lemma.
\end{abstract}
\noindent {\it Keywords and Phrases}:  Borel-Cantelli lemma; Barndorff-Nielsen lemma; limit laws.

\noindent {\it AMS 2000 Subject Classification:} 60G70, 62G30

\section{ Introduction}
Suppose $A_1,A_2,\ldots$ is a sequence of events on a common probability space and that $A^c_i$ denotes the complement of event $A_i$. The Borel-Cantelli lemma, presented below as Lemma
\ref{lemma1.1}, is used  for producing strong limit results.
\begin{lemma}\label{lemma1.1}
\begin{enumerate}
\item If, for any sequence  $A_1,A_2,\ldots$ of events,
\begin{equation}\label{1.1}
\sum_{n=1}^\infty P(A_n)<\infty,
\end{equation}
then $P(A_n\ i.o.)=0$;
 \item If $A_1,A_2,\ldots$ is a sequence
of independent events and if $\sum_{n=1}^\infty P(A_n)=\infty$,
then $P(A_n\ i.o.)=1$.
\end{enumerate}
\end{lemma}
The independence condition in the second part of Lemma~\ref{lemma1.1} has been weakened by a number of authors, including
Chung and Erdos (1952), Erdos and Renyi (1959), Lamperti (1963), Kochen and Stone (1964), Spitzer (1964), Chandra (1999, 2008),  Petrov (2002, 2004), Frolov (2012) and others.

The first part of  Borel-Cantelli lemma has been generalized in Barndorff-Nielsen (1961),  Balakrishnan and Stepanov (2010) and Frolov (2014). For
a review on the Borel–-Cantelli lemma, one may refer to the book of Chandra (2012). The result of Barndorff-Nielsen is presented below as Lemma~\ref{lemma1.2}.
\begin{lemma}\label{lemma1.2}
Let $A_n\ (n\geq 1)$ be a sequence of events such that $P(A_n)\rightarrow 0$. If
\begin{equation}\label{1.2}
\sum_{n=1}^\infty P(A_n A^c_{n+1})<\infty,
\end{equation}
then $P(A_n\ i.o.)=0$.
\end{lemma}
Observe  that condition (\ref{1.2})  in Lemma~\ref{lemma1.2} is weaker than  condition (\ref{1.1}) in the Borel-Cantelli lemma. 

In this note, we show that the  result  stated in the Barndorff-Nielsen lemma  does not follow from the argument proposed in the proof. That way, although the result is correct, the proof presented in   Barndorff-Nielsen (1961) is incomplete and not rigorous.

In connection with the above, we propose an alternative proof of this lemma based on our earlier result, presented  as Lemma~\ref{lemma2.1}. Lemma~\ref{lemma2.1}, in its turn, follows from a more general  result, given here as Lemma~\ref{lemma2.2} and proved in Balakrishnan and Stepanov (2010). In the end of this short note, we also propose a new  generalization of  Borel-Cantelli lemma.

\section{On the Proof of Barndorff-Nielsen Lemma}
First, we cite the proof of Lemma~\ref{lemma1.2} from Barndorff-Nielsen (1961) and discuss it.
\newline
\newline
\noindent {\em "Since
\begin{equation}\label{2.0}
P(A_n^c\ i.o.)=\lim_{n\rightarrow \infty}P\left(\cup_{i=n}^\infty A_i^c\right)\geq \lim_{n\rightarrow \infty}P(A_n^c)=1
\end{equation}
we have, in consequence of (\ref{1.2}) and the Borel-Cantelli lemma,
\begin{equation}\label{2.1}
P(A_n\ i.o.)=P(A_n\cap A^c_{n+1}\ i.o.)=0."
\end{equation}}
\newline
\newline
Observe that  (\ref{2.0}) implies only that
$$
P(A_n\ i.o.)=P(A_n\ i.o.\ \cap A_{n+1}^c\ i.o.).
$$
On p. 151 of  Shiryaev (1989), one can find that for some sequences of events $B_n, C_n\ (n\geq 1)$
$$
P(B_nC_n\ i.o.)< P(B_n\ i.o.\ \cap C_n\ i.o.).
$$
The last inequality shows us that  Lemma~\ref{lemma1.2} does not follow directly from the argument proposed in its proof.   We then present an alternative rigorous proof of Lemma~\ref{lemma1.2}. 
\begin{gproof}{of Lemma~\ref{lemma1.2}} 
In order to prove Lemma~\ref{lemma1.2}, we adduce  Lemma~\ref{lemma2.1}.
\begin{lemma}\label{lemma2.1}
Let $A_n,\ (n\geq 1)$ be a sequence of events such that
$P(A_n)\rightarrow 0$.  If
\begin{equation}\label{2.2}
\sum_{n=1}^\infty P(A^c_n A_{n+1})<\infty,
\end{equation}
then $P(A_n\ i.o.)=0$.
\end{lemma}
\noindent Now, the result of Lemma~\ref{lemma1.2} readily follows from Lemma~\ref{lemma2.1} and the identity
$$
\sum_{n=1}^\infty P(A_nA_{n+1}^c)=P(A_1)+\sum_{n=1}^\infty P(A^c_nA_{n+1}).
$$
\end{gproof}
Lemma~\ref{lemma2.1}, in its turn,  follows from a more general result, Lemma~\ref{lemma2.2} obtained in  Balakrishnan and Stepanov (2010). For the convenience of reader, the last lemma and its proof are also presented here.
\begin{lemma}\label{lemma2.2}
Let  $A_1,A_2,\ldots$ be a sequence of events such that $P(A_n)\rightarrow 0$.  If, for some $m\geq0$,
$$
\sum_{n=1}^\infty P(A^c_nA_{n+1}^c\ldots A_{n+m-1}^c A_{n+m})<\infty,
$$
then $P(A_n\ i.o.)=0$.
\end{lemma}
\begin{gproof}{of Lemma~\ref{lemma2.2}} 
We first note that
$$
P(A_n\ i.o.)=P\left(\bigcap_{n=1}^\infty \bigcup_{k=n}^\infty A_k\right)=\lim_{n\rightarrow\infty}P\left(\bigcup_{k=n}^\infty A_k\right).
$$
However,
\begin{eqnarray*}
P\left(\bigcup_{k=n}^\infty A_k\right)&=&P(A_n)+P(A^c_nA_{n+1})+P(A^c_nA^c_{n+1}A_{n+2})+\ldots\\
&\leq &P\{A_n\}+P(A^c_nA_{n+1}\}+\ldots+P(A^c_n\ldots A^c_{n+m-2}A_{n+m-1})\\
&+&\sum_{k=n}^\infty P(A^c_k\ldots A^c_{k+m-1}A_{k+m})\rightarrow 0
\end{eqnarray*}
as $n\rightarrow\infty$. Hence, the result follows.
\end{gproof}
Observe that Lemma~\ref{lemma2.1} follows from   Lemma~\ref{lemma2.2}, if we choose $m=1$.

\section{New Result}
In this section, we present a new theoretical result, which generalizes the second part of  Borel-Cantelli lemma. We first introduce a new notion.
\begin{definition}\label{definition3.1}
Let $A$ and $B$ be some events. We say that $\alpha\geq 0$ is the power-$A$ coefficient of dependence between $A$ and $B$ if $P(AB)=(P(A))^\alpha P(B)$, provided that $\alpha=1$ if $P(B)=0$. 
\end{definition}
\noindent Obviously, if $A$ and $B$ are independent, then $\alpha =1$. 

Let now $\bar{A}_n^c =A_n^cA_{n+1}^c\ldots$ and $\alpha _n$ be the power-$A_n^c$ coefficient of dependence between $A_n^c$ and $\bar{A}_{n+1}^c$. In the following lemma we present  sufficient conditions for $P(A_n\ i.o.)=1$.
\begin{lemma}\label{lemma3.1}
Let  $A_1,A_2,\ldots$ be a sequence of events. If
\begin{equation}\label{3.1}
\sum_{n=1}^\infty \alpha _nP(A_n)=\infty,
\end{equation}
then
\begin{equation}\label{3.2}
P(A_n\ i.o.)=1.
\end{equation}

\end{lemma}
\begin{gproof}{of Lemma~\ref{lemma3.1}} Indeed,
\begin{eqnarray*}
P(\bar{A}_n^c)&=&(P(A_n^c))^{\alpha _n}P(\bar{A}_{n+1}^c)\\
&=&\ldots =  (P(A_n^c))^{\alpha _n}\ldots(P(A_{n+k-1}^c))^{\alpha _{n+k-1}}P(\bar{A}_{n+k}^c)\quad (n, k\geq 1).
\end{eqnarray*}
By the inequality $\log (1-x)\leq  -x\ (0\leq x<1)$, we get  that
$$
P(\bar{A}_n^c)\leq  e^{-\sum_{i=n}^{n+k-1}{\alpha _iP(A_i)}}P(\bar{A}_{n+k}^c)\quad (n, k\geq 1).
$$
Then
\begin{equation}\label{3.3}
\lim_{k\rightarrow \infty }P(\bar{A}_n^c)=P(\bar{A}_n^c)\leq e^{-\sum_{i=n}^{\infty }{\alpha _iP(A_i)}}(1-P(A_n\ i.o.))\quad (n\geq 1).
\end{equation}
Obviously, (\ref{3.1}) and (\ref{3.3}) imply (\ref{3.2}).
\end{gproof}

It should be noted that Lemma~\ref{lemma3.1} is more of  a theoretical value. Indeed,  it is not easy to find the coefficients $\alpha_n$ in  general case. In some situations, however, it is not difficult. Let us consider one of such cases. Let $I_A$ be the indicator-function of event $A$, i.e. 
\begin{eqnarray*}
I_A=\left\{ \begin{array}{cc}1 , & if\ \  A\ happens,\\ 0, & otherwise.\\
\end{array}
\right.
\end{eqnarray*}
\begin{definition}\label{definition3.2}
We say that $A_n\ (n\geq 1)$ is a Markov sequence of events if the sequence of random variables $I_{A_n}\ (n\geq 1)$ is a Markov  chain. 
\end{definition}
\noindent In respect to this definition, see also Stepanov (2014).
\begin{remark}\label{remark3.1}
Let $\beta _n$ be the power-$A_n^c$ coefficient of dependence of $A_n^c$ and $A_{n+1}^c$. Then $\alpha_n=\beta_n$.
\end{remark}
In the case when  $A_n\ (n\geq 1)$ forms the Markov sequence of events, Remark~\ref{remark3.1} allows to find $\alpha _n$ easily. Indeed,
$$
\alpha _n=\frac{\log P(A_n^c\mid A_{n+1})}{\log P(A_n^c)}.
$$
\section*{References}

\begin{description} 
\item Balakrishnan, N., Stepanov, A. (2010).\ Generalization of Borel-Cantelli lemma. {\it The Mathematical Scientist}, {\bf 35}, 61--62.

\item Barndorff-Nielsen, O. (1961).\ On the rate of growth of the partial maxima of a sequence of independent identically distributed random variables. {\it Math. Scand.},  9, 383--394.

\item Chandra, T.K. (1999).\ {\it A First Course in Asymptotic Theory of Statistics}. Narosa Publishing House Pvt. Ltd., New Delhi.

\item Chandra, T.K.  (2008).\ Borel-Cantelli lemma under dependence conditions. {\it Statist. Probab. Lett.}, {\bf 78}, 390–-395.

\item Chandra, T.K., (2012).\ {\it The Borel–-Cantelli Lemma}. Springer Briefs in Statistics.

\item Chung, K.L. and Erdos, P. (1952).\ On the application of the Borel-Cantelli lemma. {\it Trans.  Amer. Math. Soc.},  72,
179--186.

\item Erdos, P. and Renyi, A. (1959).\ On Cantor's series with convergent $\sum 1/q_n$. {\it Ann. Univ. Sci. Budapest. Sect.
Math.},  2, 93--109.

\item Frolov, A.N. (2012).\ Bounds for probabilities of unions of events and the Borel--Cantelli lemma. {\it Statist. Probab. Lett.}, {\bf 82}, 2189--2197.

\item Frolov, A.N.  (2014).\ On inequalities for probabilities of unions of events and the Borel-Cantelli lemma. {\it Vestnik St. Petersb. Univ. Math.}, 47, 68–-75.

\item Kochen, S.B. and Stone, C.J. (1964).\ A note on the Borel-Cantelli lemma. {\it Illinois J. Math.},  8, 248--251.

\item Lamperti, J. (1963).\ Wiener's test and Markov chains. {\it J. Math. Anal. Appl.},  6, 58--66.

\item Petrov, V.V. (2002).\ A note on the Borel-Cantelli lemma. {\it Statist. Probab. Lett.}, {\bf 58}, 283--286.

\item  Petrov, V.V. (2004).\ A generalization of the Borel--Cantelli lemma. {\it Statist. Probab. Lett.}, {\bf 67}, 233--239.

\item Spitzer, F. (1964).\ {\it Principles of Random Walk}.  Van Nostrand, Princeton, New Jersey.

\item Shiryaev, A. (1989).\ {\it Probability}. Moscow, Nauka (in Russian).

\item Stepanov A. (2014).\ On the Use of the Borel-Cantelli Lemma   in Markov Chains, {\it Statist. Probab. Lett.}, {\bf 90}, 149--154.
\end{description}

\end{document}